\documentclass[11pt]{amsart}
\usepackage[colorlinks=true, urlcolor=blue,bookmarks=true,bookmarksopen=true,
citecolor=blue,hypertex]{hyperref}
\usepackage[curve]{xypic}
\usepackage{graphicx}

\newcounter{warn}

\newtheorem{theorem}{Theorem}[section]

\newtheorem{corollary}[theorem]{Corollary}

\newtheorem{lemma}[theorem]{Lemma}

\theoremstyle{remark}
\newtheorem{remark}{Remark}

\newcommand{\Z}{\mathbb{Z}}
\newcommand{\Q}{\mathbb{Q}}
\newcommand{\R}{\mathbb{R}}
\newcommand{\C}{\mathbb{C}}
\newcommand{\CP}{\mathbb{CP}}
\newcommand{\Ham}{\mathrm{Ham}}

\newcommand{\cmin}{c_\mathrm{min}}
\newcommand{\Ring}{\mathcal{R}}
\newcommand{\dChain}{\partial}
\newcommand{\dFloer}{\mathrm{d}}
\newcommand{\TrivLoops}{\mathcal{L}}
\newcommand{\CovTrivLoops}{\tilde{\TrivLoops}}
\newcommand{\Orbits}{\mathcal{P}}
\newcommand{\CovOrbits}{\tilde{\Orbits}}
\newcommand{\CovGrp}{\Gamma_{0}}
\newcommand{\action}{\mathfrak{a}}
\newcommand{\M}{\mathcal{M}}

\newcommand{\Mbar}{\overline{\M}}
\newcommand{\spinemap}{\sigma}

\begin{document}

\title{Quantization of the Serre spectral sequence} \author[J.-F. Barraud]{Jean-Francois Barraud} \author[O. Cornea]{Octav Cornea}
\begin{abstract}
The present paper is a continuation of \cite{BaCo1} and \cite{BaCo2}.
 It explores how the spectral sequence introduced in \cite{BaCo1}
interacts with the presence of bubbling. As consequences are obtained
 some relations between binary Gromov-Witten invariants and relative
Ganea-Hopf invariants, a criterion for detecting the monodromy of
bubbling as well as algebraic criteria for the detection of periodic
orbits.
\end{abstract}

\maketitle
\tableofcontents
\section{Introduction}
In \cite{BaCo1} has been introduced an algebraic way to encode the
properties of high dimensional moduli spaces of trajectories in
Morse-Floer type theories. The basic idea is that, by making use of a
 ``representation'' theory of the relevant moduli spaces
$$
\mathcal{M}(x,y)\stackrel{l_{x,y}}{\longrightarrow} G
$$
into some sufficiently large topological monoid $G$, one can define a
 ``rich'' Morse type chain complex whose differential is of the usual
 form
$$
dx=\sum_{y}a_{x,y}y
$$
but $a_{x,y}$, the coefficient ``measuring'' the moduli space
$\mathcal{M}(x,y)$, belongs to a graded ring (for example, the ring
of cubical chains of $G$) and is in general not zero when $\dim
(\mathcal{M}(x,y))>0$. By representation theory it is meant here not
only that the maps $l_{x,y}$ are continuous but also that they are
compatible in the obvious way with compactification and with the
crucial boundary formula:
\begin{equation}\label{eq:bdry0}
\partial \overline{\mathcal{M}}(x,y)=\bigcup_{z}\overline{\mathcal{M}}(x,z)\times
\overline{\mathcal{M}}(z,y)~.~
\end{equation}
The complex constructed this way comes with a natural filtration
induced by the grading of the generators $x,y,...$. The pages of
order greater than $1$ of the associated spectral sequence are
invariant with respect to the various choices made in the
construction and their differentials encode algebraically the
properties of the $\mathcal{M}(x,y)$'s.

This construction is described in the absence of bubbling in
\cite{BaCo1} and, in \cite{BaCo2}, it is shown to be easily
extendable to cases when pseudo-holomorphic spheres and disks exist
as long as we work under the threshold of bubbling.

\

The present paper explores what happens when bubbling does occur.

\

It is obvious that to study this case it is natural to start with the
 Hamiltonian version of Floer homology and this is indeed the setting
 of this paper. In particular, the moduli spaces $\mathcal{M}(x,y)$
consist of Floer tubes and the monoid $G$ is the space of pointed
Moore loops on $M$, $\Omega M$, with $(M^{2n},\omega)$ our underlying
 symplectic manifold.  We will also restrict to the monotone case
even  if the machinery described here appears to extend to the
general  case. The reason for this is that the main phenomena we have
identified are already present in this case and, at the same time, in
this way we avoid to deal with the well-know transversality issues
which are present in full generality.

\

Here is a short summary of our findings. Firstly, it is not
surprising that when bubbling is possible, only some of the pages of
the spectral sequence mentioned before exist. It is also expectable
that the number of pages that are defined should roughly be the
minimal Chern class, $c_{min}$, and that, moreover, some of these
pages should again be independent of the choices made in the
construction.

What is remarkable is that, in general, these pages do not
coincide with those associated to a Morse function: a quantum
deformation is generally present. Given that in the Morse case the
resulting spectral sequence is, as shown in \cite{BaCo1}, the Serre
spectral sequence of the path-loop fibration over $M$, we see that
this construction provides a new symplectic invariant which consists
of the first $c_{min}$ pages (together with their differential) of a
spectral sequence which is a quantum deformation of the Serre
spectral sequence.  One additional important point is that, on the
last defined page, the presumptive  differential, $d^{r}$, is still
defined and invariant but might not verify $(d^{r})^{2}=0$.

Of course, the next stage is to understand - at least in part - this
quantum deformation in terms of classical Gromov-Witten invariants.
In this respect we obtain that the first interesting differential
$d^{r}$, can be expressed in terms of binary Gromov-Witten invariants
 (these are those associated to spheres with two marked points) and
Ganea-Hopf invariants (these control the classical part of the
differential). Moreover, in this case, the relation $(d^{r})^{2}=0$
becomes a relation between these two types of invariants which takes
place in the Pontryaguin ring $H_{\ast}(\Omega M)$. Undoubtedly,
this is just a first step towards understanding the deeper
relationships between the combinatorics of Gromov-Witten
invariants and classical  algebraic topology invariants encoded
in the ring structure of $H_{\ast}(\Omega M)$.

The next interesting point is to understand what happens for the
first $r$ when $(d^{r})^{2}\not=0$. Clearly, the culprit is bubbling
but interestingly enough what this non-vanishing relation detects is
monodromy - the fact that in the appropriate moduli space the
attachment point of the bubbles turns non-trivially around Floer
cylinders - which turns out to interfere with the representation maps
 $l_{x,y}$. The fact that $d^{r}$ is invariant but, simultaneously,
$(d^{r})^{2}$ might not vanish is quite remarkable and, indeed, this
morphism $d^{r}$ has interest in itself and is seen to be, in fact, a
 generalization of the Seidel homomorphism \cite{Seidel}. Finally, we also discuss
an application of this structure to the detection of periodic orbits.
 This provides a sort of algebraic counterpart to the result of
Hofer-Viterbo \cite{HoVi}.

\

The paper is structured as follows. In the second section we
introduce the main notation and give the precise statements of our
results. The third section contains the proofs. In the last section we first
shortly mention some possible extensions of the construction we then
provide some examples and, finally, we  discuss the aplication to
periodic orbits.

\subsection{Acknowledgment.}
It is our great pleasure to dedicate this paper to Dusa McDuff on the 
occasion of her 60st birthday. This is even more
appropriate as, early in this project, we believed, {\em a
posteriori} without justification, that the monodromy of bubbling is
much less relevant and it is one of Dusa's questions which made us
reconsider the issue and appreciate the full importance of this
phenomenon.

\section{Notation and statement of results}

\subsection{Setting and recalls}
Fix the symplectic manifold $(M^{2n},\omega)$ and we suppose for now that $M$ is closed.
We assume that $M$ is monotone in the sense that the two morphisms $\omega:\pi_{2}(M)\to
\mathbb{R}$ and $c_{1}:\pi_{2}(M)\to \mathbb{Z}$ are proportional
with a positive constant of proportionality $\rho$. We denote by
$c_{min}$ the minimal Chern class and by $\omega_{min}$ the
corresponding minimal symplectic area (so that we have
$\omega_{min}=\rho c_{min}$).

\subsubsection{Binary Gromov-Witten invariants.}
Fix on $M$ a generic almost complex structure $J$ which tames
$\omega$. The binary Gromov Witten invariants we are interested in
can be described as follows: pick a generic Morse function $f$ and
metric on $M$. Denote by $i(x)=ind_{f}(x)$ for each $x\in Crit(f)$.
For two critical points $x$ and $y$ and a class $\alpha\in
\pi_{2}(M)$ such that $i(x)-i(y)+2c_{1}(\alpha)-2=0$, we define
$GW_{\alpha}(x,y)$ as the number of elements in the moduli space
$\mathcal{M}(J,\alpha;x,y)$ which consists of $J$-holomorphic spheres
 in the homology class $\alpha$ with two marked points, one lying on
the unstable manifold of $x$ and the other on the stable manifold of
$y$, modulo reparametrization. As such $GW_{\alpha}(x,y)$ is not an
invariant (because $x,y$ might not be Morse cycles). However, if for
two classes $[x]=[\sum\lambda_{i}x_{i}]$ and $[y]=[\sum \mu_{i}
y_{i}]$, we define
$GW_{\alpha}([x],[y])=\sum\lambda_{i}\mu_{j}GW_{\alpha}(x_{i},y_{j})$
 then we obtain an invariant. For $\alpha\in\pi_{2}$, let $[\alpha]$
be its image by the morphism $\pi_{2}(M)\to H_{1}(\Omega M)$.

\subsubsection{The Novikov ring.}
Let $\TrivLoops(M)$ be the space of contractible loops in $M$.

Let $\Gamma$ be the image of the Hurevicz morphism $\pi_{2}(M)\to
H_{2}(M,\Z/2)$. The two forms $\omega$ and $c_{1}$ define morphisms
$\Gamma:\xrightarrow{\omega , c_{1}}\R , \Z$ which under our
monotonicity assumption are proportional. Let $\CovGrp
=\Gamma/\ker(\omega)$. We let $\Lambda$ be the associated Novikov
ring which is defined as follows
$$
  \Lambda=\left\{             \sum_{\alpha\in\CovGrp
  }\lambda_{\alpha}e^{\alpha}\           \right\}.
$$
where the coefficients $\lambda_{\alpha}$ belong to $\Z/2$ such that
$$
\forall c>0,\ \sharp\{\alpha,\lambda_{\alpha}\neq0,\omega(\alpha)\leq
c\}<+\infty~.~
$$
The grading of the elements in $\Lambda$ is given by
$|e^{\lambda}|=-2c_{1}(\lambda)$.

We also denote by $\CovTrivLoops(M)$ the covering of $\TrivLoops(M)$
associated to $\Gamma_{0}$: it is the quotient of the space of
couples $(\gamma,\Delta)$, where $\gamma\in\TrivLoops(M)$ and
$\Delta$ is a disk bounded by $\gamma$, under the equivalence
relation $(\gamma,\Delta)\sim(\gamma',\Delta')$ if $\gamma=\gamma'$
and $\omega([\Delta-\Delta'])=c_{1}([\Delta-\Delta'])=0$.

\begin{remark}
Here and later in the paper we could also use, alternatively,
rational coefficients as all the moduli spaces involved are
orientable and the orientations are compatible with our formulae.
\end{remark}

\subsubsection{Moduli spaces of Floer tubes.}\label{subsubsec:moduli_Floer_tub}

Let $H:M\times S^{1}\to\R$ be a Hamiltonian function. The Hamiltonian
flow associated to $H$ is the flow of the (time dependent) vector
field $X_{H}$ defined by~:
$$
\omega(X_{H_{t}},\cdot)= -dH_{t}
~.~$$
All along this paper, the periodic orbits of $X_{H}$ will be supposed
 to be non degenerate. We denote by $\Orbits_{H}\subset\TrivLoops(M)$
 the set of all contractible periodic orbits of the hamiltonian flow
associated to $H$ and we let $\CovOrbits_{H}$ be the covering of
$\Orbits_{H}$ which is induced from $\CovTrivLoops(M)$.

For each periodic orbit $x\in \Orbits_{H}$ we fix a lift
$(x,\Delta_{x})\in \CovOrbits_{H}$. For a generic pair $(H,J)$ and
$x,y\in \Orbits_{H}$, $\lambda\in \CovGrp$ we now consider the moduli
spaces:
$$
\mathcal{M}'(x,y;\lambda)=\{u:\R\times S^{1} : u\ \textrm{verifies} \
 (\ref{eq:floer})\}
$$
so that the pasted sphere $\Delta_{x}\cup u\cup (-\Delta_{y})$ is of
class $\lambda$ and
\begin{equation}\label{eq:floer}
 \partial_{s} u + J(u)\partial_{t}u-J(u)X_{H}(u)=0\ , \lim_{s\to -\infty}u(s,t)=x(t), \lim_{s\to+\infty}u(s,t)=y(t)~.~
\end{equation}
Of course, these moduli spaces are quite well-known in the subject
and we refer to \cite{Sa} for their properties. In
particular, they have natural orientations and, when
$(x,\Delta_{x})\not=(y,\Delta_{y})$ they admit a free $\R$ action. We
denote the quotient by this action by $\mathcal{M}(x,y;\lambda)$ and
we have
$$
\dim \mathcal{M}(x,y;\lambda)=\mu((x,\Delta_{x}))-\mu((y,\Delta_{y}))+2c_{1}(\lambda)-1 %
$$
where $\mu((x,\Delta_{x})$ is the Conley-Zehnder index of the orbit
$x$ computed with respect to the capping disk $\Delta_{x}$.

\subsubsection{Monodromy of bubbling}
Among the standard properties of the moduli spaces above we recall
that they admit a natural topology as well as natural
compactifications, $\overline{\mathcal{M}}(x,y;\lambda)$, so that the
 following formula is valid:

\begin{equation}\label{eq:bdry}
\partial \overline{\mathcal{M}}(x,y;\lambda)=\bigcup_{z,\lambda'+\lambda''=
\lambda}\overline{\mathcal{M}}(x,z;\lambda')\times
\overline{\mathcal{M}}(z,y;\lambda'')\ \cup\ \Sigma_{x,y,\lambda}
\end{equation}

Here $\Sigma_{x,y,\lambda}$ is a set of codimension $2$ which
consists of Floer tubes with at least one attached bubble.

\

We will say that $(H,J)$ has {\em bubbling monodromy} if there exist
$x,y\in \Orbits_{H}$ and $\lambda\in\CovGrp$ so that:
$$
H^{1}(\Sigma_{x,y,\lambda};\Z)\not=0
$$
This means, in particular, that $\pi_{1}(\Sigma_{x,y,\lambda})\not=0$
 so that there are non-contractible loops in the space of Floer tubes
 with bubbles.

\subsubsection{Truncated differentials and spectral sequences.}
\label{subsubsec:algebra_truncated} The following algebraic notions
will be useful in the formulation of  our results.

We say that the sequence of graded vector spaces $(E^{r},d^{r})$,
$0\leq r\leq k$ is a {\em truncated spectral sequence of order $k$}
if $(E^{r}, d^{r})$ is a chain complex for each $r\leq k-1$ which
verifies $H_{\ast}(E^{r},d^{r})=E^{r+1}$ and $d^{k}$ is a linear map
of degree $-1$. A truncated spectral sequence of $\infty$-order is a
usual spectral sequence. A morphism of order $k$ truncated spectral
sequences is a sequence of chain maps $\phi_{r}:(E^{r},d^{r})\to
(F^{r},d^{r})$, $0\leq r\leq k$, so that
$H_{\ast}(\phi_{r})=\phi_{r+1}$ for $0\leq r\leq k-1$. We say that
two truncated spectral sequences are isomorphic starting from page
$s$ is they admit a morphism which is an isomorphism on page $s$
(and, hence, on each later page).

The typical example of a truncated spectral sequence appears as
follows. Assume that $C_{\ast}$ is a graded rational vector space and
 that $F^{i}C$ is an increasing filtration of $C_{\ast}$. We say that
 a linear map $d:C_{\ast}\to C_{\ast-1}$ is a {\em truncated
differential of order $k$} compatible with the given filtration if
$d(F^{i}C)\subset F^{i}C$ $\forall i$ and
$$
(d \circ d)(F^{r}C) \subset F^{r-2k}C
$$ for all $r\in \Z$.
It is easy to see that a truncated differential of order $k$ induces
a truncated spectral sequence of the same order. Indeed, by using the
standard descriptions of the $r$-cycles
$$
Z^{r}_{p}=\{ v\in F^{p}C : dv\in F^{p-r}C\}+F^{p-1}C
$$
and $r$-boundaries
$$
B^{r}_{p}=\{d F^{p+r-1}C \cap F^{p}C\}+F^{p-1}C
$$
it is immediate to see that $B^{r}_{p}\hookrightarrow Z^{r}_{p}$ for
$0\leq r\leq k$ which allows us to define the pages of the truncated
spectral sequence by $E^{r}_{p}=Z^{r}_{p}/B^{r}_{p}$. Obviously, $d$
induces differentials $d^{r}$ on $E^{r}$ when $r< k$ as well as a
degree $-1$ linear map $d^{k}$ on $E^{k}$.

\subsection{Main statement.}
We will formulate our main statement in a simple case and we will
discuss various extensions at the end of the paper. Therefore, we
assume here that $(M,\omega)$ is closed, simply-connected and
monotone with $c_{min}\geq 2$.

\begin{theorem}
\label{theo:spec}  There exists a truncated spectral sequence of
order $c_{min}$,   $(E^{r}(M), \dFloer^{r})$, whose isomorphism type
starting from page   $2$ is a symplectic invariant of $(M,\omega)$
and which has the   following additional properties:
\begin{itemize}
\item[i.] As a bi-graded vector space we have an isomorphism:
  $$
  E^{2}\cong H_{\ast}(M)\otimes H_{\ast}(\Omega M)\otimes \Lambda~.~
  $$
\item[ii.] The differential $\dFloer^{2}$ has the decomposition
  $$
  \dFloer^{2}=\dFloer^{2}_{0}+\dFloer^{2}_{Q}
  $$
  where $\dFloer^2_0$ is the differential appearing in the classical
    Serre spectral sequence of the path loop fibration $\Omega M\to
   PM\to M$ and
  $$
  \dFloer^{2}_{Q}x=\sum_{y,\alpha}
  GW_{\alpha}(x,y)y[\alpha]e^{\alpha}~.~
  $$
\item[iii.] If $(\dFloer^{\cmin})^{2}\not=0$, then any regular pair $(H,J)$   has bubbling monodromy.
\end{itemize}
\end{theorem}
\begin{remark}
{\rm Clearly, if $\dFloer^{2}\circ \dFloer^{2}=0$ - for example if
$c_{min}\geq 3$ - the vanishing of the square of $d^{2}\circ d^{2}$
translates into some relations between binary Gromov-Witten
invariants and the classical Serre spectral sequence differential
$\dFloer_{0}^{2}$. In turn, this differential is quite well known and
 rather easy to compute and it can be expressed in many cases in
terms  of relative Ganea-Hopf invariants (see \cite{Co1}). The
interesting  part about these relations is that they take place in
the Pontryaguin  algebra $H_{\ast}(\Omega M)$. Indeed, in the formula
at ii.  $\alpha\in H_{1}(\Omega M)$, $e^{\alpha}\in \Lambda$ and
$x,y\in  H_{\ast}(M)$ so that in the square of the differential
appears the  Pontryaguin product $H_{1}(\Omega M)\otimes H_{1}(\Omega
M)\to  H_{2}(\Omega M)$.
}\end{remark}

The relation with the Seidel homomorphism is seen by considering the
spectral sequence in the case of a symplectic fibration over
$\CP^{1}$.

We also formulate here a very simple version of our application to
the detection of periodic orbits.  We specialize to the case when the
manifold $M$ admits a perfect Morse function (that is a Morse
function whose associated Morse complex has trivial differential). We
also need the following notion. Let $x,y\in H_{\ast}(M)$ and
$\lambda\in\Lambda$. We will say that $x$ and $ye^{\lambda}$ (which
exist on the $E^{2}$ page of the spectral sequence in Theorem
\ref{theo:spec}) are $\dFloer^{r}$-related if $x$ survives to the
$r$-th level of the spectral sequence and there is some $\gamma\in
C_{\ast}(\Omega M)$ so that  the product $\gamma\otimes y
e^{\lambda}$ also survives to the $r$-th page of the spectral
sequence and we have $\dFloer^{r}([x])=[\gamma\otimes y
e^{\lambda}]+\ldots$

\begin{corollary}\label{thm:extra periodic orbits}
Assume that there are homology classes $x,z\in H_{\ast}(M)$, $|x|< |z|$,
so that $x$ is $\dFloer^{r}$-related to $ze^{\lambda}$ and   $H_{k}(M)\otimes \Lambda_{q}=0$ for $|x|> k+q > |ze^{\lambda}|)$.
Then any self-indexed perfect Morse function on $M$ has some non-trivial closed characteristic.
\end{corollary}

By a self-indexed Morse function $f$ we mean here that the critical
points of the same index have the same critical value and $ind_{f}(x)> ind_{f}(y)$
implies $f(x)>f(y)$.

There are many ways in which this statement can be extended and some will be discussed
at the end of the paper.

\section{Proof of the main theorem. }

\subsection{Construction of the truncated spectral sequence.}\label{subsec:definition_ss}

In this section we fix the $1$-periodic hamiltonian $H$ and almost
complex structure $J$ compatible with $\omega$ so that the pair
$(H,J)$ is generic (of course, both are in general time-dependent).
For simplicity, we will also assume to start that the manifold $M$ is
 simply-connected but we will see later on that this condition can be
 dropped with the price that the construction becomes more complicated.

As in \cite{BaCo1} the truncated spectral sequence we intend to
discuss is induced by a natural filtration of an enriched Floer type
pseudo-complex. We use the term pseudo-complex here to mean that we
will not have here a true differential but rather a truncated one.
The construction of this pseudo-complex  is a refinement of the
classical Floer construction in which the coefficient ring is
replaced with the ring of cubical chains over the Moore loops on $M$.
Here is this construction in more detail.

\subsubsection{Coefficient rings}\label{plan:Coefficient ring}

Let $C_{*}$ denote the ``cubical'' chain complex, let $\Omega X$ be
the Moore loop space over $X$ (the space of loops parametrized by
intervals of arbitrary length). Consider the space $M'$ obtained from
 $M$ by collapsing to a point a simple path $\gamma$ going through
the  starting point of each periodic orbit. Notice that
$C_{\ast}(\Omega  M')$ is a differential ring where the product is induced by the
concatenation of  loops. Finally, our coefficient ring is:
$$\Ring_{\ast} = C_{\ast}(\Omega M')\otimes \Lambda~.~$$
This is a (non abelian) differential ring, and its differential will
be denoted by $\dChain$.

The (pseudo)- complex we are interested in is a (left) differential module
generated by the contractible periodic orbits of $H$ over this ring:
$$
C(H,J)=\oplus_{\tilde{x}\in\CovOrbits_{H}} \Ring_{*}\ \tilde{x}/\sim
$$
with the identification $\tilde{x}e^{\lambda}\sim \tilde{x}\sharp
\lambda$, where $\tilde{x}\sharp \lambda$ stands for the capping of
$x$ obtained by gluing a sphere in the class $\lambda$ to
$\tilde{x}$. The grading of an element in $\tilde{x}\in
\CovOrbits_{H}$ is given by the respective Conley-Zehnder index.
There is a natural filtration of this complex which is given by
$$
F^{r}C(H,J)=\Ring_{\ast} < \tilde{x}\in \CovOrbits_{H} : \
\mu(\tilde{x})\leq r>~.~
$$
We will call this the canonical filtration of $C(H,J)$.

\subsubsection{Truncated boundary operator}\label{plan:Boundary operator}
The next step is to introduce a truncated differential on $C(H,J)$.
We recall from \S \ref{subsubsec:moduli_Floer_tub} the definition of
the moduli spaces $\mathcal{M}(x,y;\lambda)$ of Floer tubes. We
recall also that this definition requires a choice of lift
$\tilde{x}\in \CovOrbits_{H}$ for each $x\in\Orbits_{H}$. With these
conventions and - as assumed before - for a generic choice of $J$ and
 $H$ - the moduli spaces are smooth manifolds of dimension
$|\tilde{x}|-|\tilde{y}|-1$ when $|\tilde{x}|\neq|\tilde{y}|$, and
they have a natural compactification involving ``breaks'' of the
tubes on intermediate orbits, or bubbling off of holomorphic spheres.
 We will write $\mathcal{M}(\tilde{x},\tilde{y})$ for the moduli
 space of Floer tubes which lift to paths inside $\CovTrivLoops (M)$
joining $\tilde{x}\in\CovOrbits_{H}$ to $\tilde{y}\in \CovOrbits_{H}$
 and we let $\overline{\mathcal{M}}(\tilde{x},\tilde{y})$ be the
respective compactification. In our monotone situation these
compactifications are pseudo-cycles with boundary.

\

To define the truncated boundary operator we proceed as in the usual
Floer complex, but we intend to take into consideration the moduli
spaces of arbitrary dimensions instead of restricting to the $0$
dimensional ones. To associate to the (compactification of the)
moduli spaces coefficients in our ring $\Ring$, we first need to
represent them into the loop space $\Omega (M')$, and then pick
chains representing them (i.e. defining their fundamental classes
relative to their boundary).

Let us start with ``interior'' trajectories, i.e. elements
$v\in\M(\tilde{x},\tilde{y})$. Let $u:\R\times S^{1}\to M$ be a
parametrization of $v$. Since the value of the action functional
$$
\action_{H} : \CovTrivLoops (M) \to \R,\  \action_{H}((\gamma,\Delta))=-\int_{D^{2}} \Delta^{\ast}\omega + \int_{S^{1}}H(t, \gamma(t))dt %
$$
is strictly decreasing along the $\R$ direction, it can be used to
reparametrize $u$ by the domain
$[-\action(\tilde{x}),-\action(\tilde{y})]\times S^{1}$, and the
restriction of $u$ to the interval
$[-\action(\tilde{x}),-\action(\tilde{y})]\times\{0\}$ defines a
Moore loop in $M'$. This defines a map
\begin{equation}
  \label{eq: Restriction map}
  \spinemap_{\tilde{x},\tilde{y}}:\M(\tilde{x},\tilde{y})\to\Omega(M')
\end{equation}
which is continuous. We will call it the ``spine'' map.

This map should then be extended to the compactification
$\Mbar(\tilde{x},\tilde{y})$ of $\M(\tilde{x},\tilde{y})$.

It is well-known that the objects in
$\overline{\mathcal{M}}(\tilde{x},\tilde{y})$ are constituted by
Floer trajectories possibly broken on some intermediate periodic
orbits to which might be attached some $J$-holomorphic spheres that
have bubbled off.

It is easy to see that the map $\spinemap_{\tilde{x},\tilde{y}}$
extends continuously over the part of this set where no spheres are
attached to some tube in a point belonging to the line
$\R\times\{0\}$. Indeed, as in \cite{BaCo1}, except for these types
of elements, the spine map is compatible with the breaking of Floer
tubes in the sense that the loop associated to a broken trajectory is
 the product of the loops associated to each ``tube'' component.

Let $\alpha_{min}\in \Gamma_{0}$ be the class so that
$c_{1}(\alpha_{min})= c_{min}$ (by our monotonicity assumption there
is a single such class). By using again the monotonicity assumption
we see that bubbling off  of a sphere in class $\alpha\in \Gamma_{0}$
can occur in a moduli space $\Mbar(\tilde{x},\tilde{y})$ with
$\tilde{y}\not=\tilde{x}\sharp\alpha$ only if
$$|\tilde{x}|-|\tilde{y}|\geq 2 c_{1}(\alpha) +1~.~$$

It is also important to note that bubbling of a sphere in the class
$\alpha$ is also possible inside the space $\Mbar
(\tilde{x},\tilde{x} \sharp \alpha)$. In all cases, bubbling of an
$\alpha$ sphere is never possible if $|\tilde{x}|-|\tilde{y}|\leq 2
c_{1}(\alpha)-1$.

We summarize this discussion:
\begin{lemma}\label{lem:bubbling_cond}
  The spine map $\spinemap_{\tilde{x},\tilde{y}}$ extends
  continuously to $\overline{\mathcal{M}}(\tilde{x},\tilde{y})$ if
  $$
   |\tilde{x}|-|\tilde{y}|\leq 2 \cmin - 1~.~
  $$
  In case $|\tilde{x}|-|\tilde{y}|= 2 \cmin$ and if $\sigma$ does not
     have such a continuous extension to
  $\overline{\mathcal{M}}(\tilde{x},\tilde{y})$, then
  $\tilde{y}=\tilde{x}\sharp\alpha_{min}$.
\end{lemma}

The spine map obtained in this way satisfies also a compatibility
condition which we now make explicit. If
$\mathcal{M}(\tilde{x},\tilde{z})\times
\mathcal{M}(\tilde{z},\tilde{y})\subset
\overline{\mathcal{M}}(\tilde{x},\tilde{y})$, then the restriction of
$\spinemap_{\tilde{x},\tilde{y}}$ on the set on the left of the
inclusion equals $m\circ(\spinemap_{\tilde{x},\tilde{z}}\times
\spinemap_{\tilde{z},\tilde{y}})$ where
$$
  m:\Omega M'\times \Omega M'\to \Omega M'
$$
is loop concatenation.

For pairs $(\tilde{x},\tilde{y})$ with $|\tilde{x}|-|\tilde{y}|\leq
2\cmin -1$, we use the map $\spinemap_{\tilde{x},\tilde{y}}$ to
represent the moduli spaces $\Mbar(\tilde{x},\tilde{y})$ inside the
loop space $\Omega (M')$. We then choose ``chain representatives''
$m(\tilde x,\tilde y)\in C_{\ast}(\Omega M')$, i.e. chains
generating the fundamental class of $\spinemap(\Mbar(\tilde x,\tilde
y))$ relative to its boundary, in such a way that:
\begin{equation}\label{eq:bdry_chain}
 \dChain m(\tilde x,\tilde y)=
 \sum_{|y|<|z|<|x|} m(\tilde x,\tilde z)*m(\tilde z,\tilde y)
\end{equation}
where $\ast$ is the operation induced on $C_{\ast}(\Omega M')$ by
the concatenation of loops.

The key point regarding this formula is that, under our assumption
$|\tilde{x}|-|\tilde{y}|\leq 2\cmin -1$, the compactified moduli
space $\overline{\mathcal{M}}(\tilde{x},\tilde{y})$ is a manifold
with boundary. Moreover, its boundary verifies the usual formula
valid in the absence of bubbling so that the construction of the
$m(-,-)$'s is the same as that in the non-bubbling setting. We refer
to \cite{BaCo1} for a complete discussion of  this construction.

We now define the boundary operator $\dFloer$ by:
\begin{equation}
  \label{eq:Boundary operator}
  \dFloer \tilde{x}=\sum_{1\leq |\tilde{x}|-|\tilde{y}|\leq2\cmin -1}
                       m(\tilde{x},\tilde{y})\ \tilde{y}
\end{equation}
and extend it to the full complex using the Leibnitz rule.

It is easy to check that $\dFloer$ has degree $-1$ with respect to the
 total grading and that it is compatible with the canonical
filtration. Notice first that if $\gamma\otimes \tilde{x}\in C_{\ast}(\Omega M)\otimes \tilde{\mathcal{P}}_{H}$
we have $\dFloer\circ\dFloer (\gamma\otimes\tilde{x})=(\gamma\otimes(\dFloer\circ\dFloer)(\tilde x))$. We now
compute:
\begin{align*}
  \dFloer\circ \dFloer(\tilde{x})
    &=\sum_{|\tilde{x}|-|\tilde{y}|\leq 2\cmin-1}
         \dFloer(m(\tilde{x},\tilde{y})\ \tilde{y})\\
    &=\sum_{1\leq|\tilde{x}|-|\tilde{y}|\leq  2\cmin-1}
         \dChain m(\tilde{x},\tilde{y})\ \tilde{y} +
         m(\tilde{x},\tilde{y})\ \dFloer \tilde{y}\\
    &=\sum_{\substack{
              1\leq|\tilde{x}|-|\tilde{y}|\leq 2\cmin-1\\
              |\tilde{y}|+1\leq|\tilde{z}|\leq |\tilde{x}|-1
              }}
         m(\tilde{x},\tilde{z})m(\tilde{z},\tilde{y})\ \tilde{y} +
      \sum_{\substack{
              1\leq|\tilde{x}|-|\tilde{y}|\leq  2\cmin-1\\
              1\leq|\tilde{y}|-|\tilde{z}|\leq  2\cmin-1
            }}
         m(\tilde{x},\tilde{y})m(\tilde{y},\tilde{z})\  \tilde{z}\\
    &=\sum_{\substack{
              1\leq|\tilde{x}|-|\tilde{y}|\leq  2\cmin-1\\
              |\tilde{y}|- 2\cmin+1\leq|\tilde{z}|\leq |\tilde{x}|- 2\cmin
            }}
         m(\tilde{x},\tilde{y})m(\tilde{y},\tilde{z})\  \tilde{z}
\end{align*}
and we see that $\dFloer^{2}$ drops the filtration index by at least
$ 2\cmin$. In the algebraic terms of \S
\ref{subsubsec:algebra_truncated} we obtain:

\begin{lemma}\label{lem:e2}
With the definition above, $\dFloer$ is a truncated differential of
order $\cmin$ with respect to the canonical filtration on $C(H,J)$
and thus it induces a truncated spectral sequence $E^{r}(H,J)$ of the
 same order so that
$$
E^{2}(H,J)\cong H_{\ast}(M)\otimes H_{\ast}(\Omega M)\otimes \Lambda
$$
\end{lemma}

The isomorphism in the lemma is obvious because
$E^{1}(H,J)\cong CF_{\ast}(H,J)\otimes H_{\ast}(\Omega M)$ and as $\dFloer^{1}$ only
involves the $0$ dimensional moduli spaces of Floer tubes we obtain that
$\dFloer^{1}$ is just: $d_{F}\otimes id$ where $(CF_{\ast}(H,J),
d_{F})$ is the usual Floer complex (with coefficients in the Novikov
ring $\Lambda$). Thus we have constructed our truncated spectral
sequence and have proved property i. in Theorem \ref{theo:spec}

\begin{remark}\label{rem:non-monotone}
Without the monotonicity assumption, but still assuming that the
moduli spaces in question are regular, there is no way to avoid the
bubbling phenomenon, even on low-dimensional moduli spaces. However, on
$2$-dimensional moduli spaces, the bubbling component is $0$
dimensional, and hence consists in isolated points: for each of them,
 the real line $\R\times\{0\}$ can actually be deformed to avoid the
 point where the bubble is attached. Interpolating between these
 perturbed real lines with the standard one in small neighborhood of
 the ``bubbled'' trajectories defines a spine map for $2$ dimensional
 moduli spaces, with the desired continuity and compatibility
 conditions.
\end{remark}

\subsection{Invariance of the truncated spectral sequence.}

To show invariance we will proceed along Floer's original proof by
first constructing a comparison morphism between the spectral
sequences associated to two different sets of generic data
$(H_{i},J_{i})_{i=0,1}$. We will describe the construction of this
morphism in more detail below but we only mention here one remarkable
 fact: despite the fact that in our spectral sequences we might have
$d^{\cmin}\circ d^{\cmin}\not=0$ it is still true that the morphism
$d^{\cmin}$ is invariant.

The construction uses a homotopy between $(H_{0}, J_{0})$ and
$(H_{1},J_{1})$. As in the usual Floer case, we consider a generic
homotopy between them, $(G,\bar{J})$, and, for
$\tilde{x}\in\CovOrbits_{H_{0}}$ and
$\tilde{x}'\in\CovOrbits_{H_{1}}$, we consider the moduli spaces
$\mathcal{N}(\tilde{x},\tilde{x}')$ of tubes $v:\R\times S^{1}\to M$
which lift in $\CovTrivLoops (M)$ to a path joining $\tilde{x}$ to
$\tilde{x}'$ and verify the equation:
\begin{equation}
  \label{eq:dbar+homotopy}
  \partial_{s}u+\bar{J}(s,u(s,t))(\partial_{t}u-X_{G}(s,u(s,t)))=\nu_{s}(u).
\end{equation}
The moduli space $\mathcal{N}(\tilde{x},\tilde{x}')$ has properties
similar to those of $\mathcal{M}'(-,-)$ except that it has no
$\R$-invariance. Its dimension is $|\tilde{x}|-|\tilde{x}'|$.
Clearly, bubbling of an $\alpha$-sphere inside such a moduli space is
 not possible if $|\tilde{x}|-|\tilde{x}'|\leq 2c_{1}(\alpha)-1$. As
in \S\ref{subsec:definition_ss} , sphere bubbling is the only
obstruction to extend the spine map. Assuming that
$|\tilde{x}|-|\tilde{x}'|\leq 2\cmin -1$ the spine map can therefore
be extended over these spaces in a way compatible with the spine maps
 of $(H_{0},J_{0})$ and $(H_{1},J_{1})$ (as in \cite{BaCo1}).

The chain morphism between the two (truncated)-complexes is defined
by a formula similar to \eqref{eq:Boundary operator}:
$$
\Theta(\tilde{x})=\sum_{0\leq|\tilde{x}|-|\tilde{x}'|\leq 2\cmin-1}
                   m'(\tilde{x},\tilde{x}') \tilde{x}'
$$
where $m'(\tilde{x},\tilde{x}')$ is a chain in the loop space
representing the moduli space $\mathcal{N}(\tilde{x},\tilde{x}')$ (as
 in \cite{BaCo1}). This morphism clearly respects the canonical
filtrations.

We also have:
$$
 \dChain m'(\tilde{x},\tilde{x}')=
 \sum_{|\tilde{x}'|\leq|\tilde{y} |\leq|\tilde{x}|-1}        m
 (\tilde{x},\tilde{y })m'(\tilde{y} ,\tilde{x}')
 +\sum_{|\tilde{x}'|+1\leq|\tilde{y}'|\leq|\tilde{x}|}
 m'(\tilde{x},\tilde{y}')m (\tilde{y}',\tilde{x}')
$$

Computing $\dFloer\Theta$ and $\Theta\dFloer$ we get:
\begin{align*}
 \dFloer\Theta(\tilde{x})
  &=\dFloer\Big(\sum_{
          0\leq|\tilde{x}|-|\tilde{x}'|\leq 2\cmin-1}
       m'(\tilde{x},\tilde{x}')\tilde{x}'\Big)\\
  &=\sum_{\substack{
          0\leq |\tilde{x} |-|\tilde{x}'|  \leq 2\cmin-1\\
          1\leq |\tilde{x} |-|\tilde{y} |\leq 2\cmin-1\\
          0\leq |\tilde{y} |-|\tilde{x}'|  \leq 2\cmin-1}}
       m (\tilde{x},\tilde{y })m'(\tilde{y} ,\tilde{x}')\tilde{x}'
   +\\&\quad
   +\sum_{\substack{
          0\leq|\tilde{x} |-|\tilde{x}'|  \leq 2\cmin-1\\
          0\leq|\tilde{x} |-|\tilde{y}'|  \leq 2\cmin-1\\
          1\leq|\tilde{y}'|-|\tilde{x}'|\leq 2\cmin-1}}
       m'(\tilde{x},\tilde{y'})m (\tilde{y}',\tilde{x}')\tilde{x}'
   +\\&\quad
   +\sum_{\substack{
          0\leq|\tilde{x} |-|\tilde{x}'|  \leq 2\cmin-1\\
          1\leq|\tilde{x}'|-|\tilde{y}'|\leq 2\cmin-1}}
       m'(\tilde{x},\tilde{x}')m(\tilde{x}',\tilde{y}')\tilde{y}'
\intertext{and}
 \Theta\dFloer(\tilde{x})
  &=\sum_{\substack{
          1\leq|\tilde{x}|-|\tilde{y} | \leq 2\cmin-1\\
          0\leq|\tilde{y}|-|\tilde{x}'|   \leq 2\cmin-1}}
       m(\tilde{x},\tilde{y})m'(\tilde{y},\tilde{x}') \tilde{x}'
\intertext{so that}
  \dFloer\Theta-\Theta\dFloer
  &=\sum_{\substack{
          0\leq|\tilde{x} |-|\tilde{x}'|  \leq 2\cmin-1\\
          |\tilde{x}'|-|\tilde{y}'|\leq 2\cmin-1\\
          |\tilde{x}|-|\tilde{y}'|\geq 2\cmin}}
       m'(\tilde{x},\tilde{x}')m(\tilde{x}',\tilde{y}')\tilde{y}'
 \\&  -\sum_{\substack{
          0\leq|\tilde{x}|-|\tilde{y} |-1 \leq 2\cmin-1\\
          |\tilde{y} |-|\tilde{x}'|\leq 2\cmin-1\\
          |\tilde{x} |-|\tilde{x}'|\geq 2\cmin}}
       m(\tilde{x},\tilde{y})m'(\tilde{y},\tilde{x}') \tilde{x}'
\end{align*}
which is not $0$, but has degree at least $-2\cmin$ with respect to
the Maslov index.

\medskip

It is easy to see that this implies that $\Theta$ induces a morphism
of truncated spectral sequences:
$$\bar\Theta : E(H_{0},J_{0})\to E(H_{1},J_{1})~.~$$

Similarly to the isomorphism in Lemma \ref{lem:e2} it is easy to see
that $E^{1}(\Theta)$ is identified with:
$$
\theta_{F}\otimes id : CF(H_{0},J_{0})\otimes H_{\ast}(\Omega M)\to
CF(H_{1},J_{1})\otimes H_{\ast}(\Omega M)
$$
where $\theta_{F}$ is the Floer comparison morphism. As this morphism
 induces an isomorphism in homology we deduce that $E^{2}(\Theta)$
and  hence all of $\bar\Theta$ are isomorphisms for $r\geq 2$ and
this  shows the invariance claim in the statement of Theorem
\ref{theo:spec}.

\begin{remark}
A morphism of spectral sequences preserves the bi-degree, therefore
to show that $\bar{\Theta}$ is a morphism we only need that
$\dFloer\Theta-\Theta\dFloer$ drops the filtration degree by $\cmin$.
 In other words, a considerable part of the geometric information
carried by $\Theta$ is actually forgotten in the spectral sequence.
There are some ways to recover it but as this goes beyond the purpose
 of the present paper we will not discuss this here.
\end{remark}

\subsection{Detection of monodromy.}\label{subsec:monodromy}
The purpose here is to prove Theorem \ref{theo:spec} iii. thus we fix
 a regular pair $(H,J)$ and we assume that $\dFloer^{\cmin}\circ
\dFloer^{\cmin}\not=0$.

We start by looking again at the calculation for $\dFloer\circ
\dFloer$ given before Lemma \ref{lem:e2}. We see from that formula
that $\dFloer^{\cmin}\circ \dFloer^{\cmin}$ is given by a linear
combination of terms of the form
$$
S(\tilde{x})=\sum_{\substack{
              0\leq|\tilde{x}|-|\tilde{y}|\leq  2\cmin-1\\
              |\tilde{y}|- 2\cmin+1\leq|\tilde{z}|=|\tilde{x}|-
              2\cmin             }}
              m(\tilde{x},\tilde{y})m(\tilde{y},\tilde{z})\
              \tilde{z}
$$

For each fixed $\tilde{z}$ with $|\tilde{z}|=|\tilde{x}|-2\cmin$ this
 last sum can be rewritten as
$$
  S(\tilde{x})=\sum_{\tilde{z}}S(\tilde{x},\tilde{z})
$$
with
$$
  S(\tilde{x},\tilde{z})=\sum_{|\tilde{x}|-1\geq |\tilde{y}|\geq | \tilde{x}|-2\cmin -1}m(\tilde{x},\tilde{y})m(\tilde{y},\tilde{z})\ \tilde{z}~.~%
$$

Suppose that $\tilde{z}\not=\tilde{x}\sharp \alpha_{min}$. In that
case, as indicated in Lemma \ref{lem:bubbling_cond}, the spine map is
 well defined and continuous on the whole space
$\overline{\mathcal{M}}(\tilde{x},\tilde{z})$ and no bubbling is
possible inside this space. But this means that we may find a
representing chain $m(\tilde{x},\tilde{z})$ so that, as in formula
(\ref{eq:bdry_chain}),
$$
(\partial m(\tilde{x},\tilde{z}))\ \tilde{z} = S(\tilde{x},\tilde{z})
$$
which means that $S(\tilde{x},\tilde{z})$ vanishes in $E^{r}$ for
$r\geq 2$.

Thus, the only terms which count in $d^{\cmin}\circ d^{\cmin}$ are
$S(\tilde{x},\tilde{x}\sharp \alpha_{min})$ and if $d^{\cmin}\circ
d^{\cmin}\not=0$, then at least one such term survives to $E^{r}$. To
 simplify notation we let $\tilde{x}\sharp \alpha_{min}=
\tilde{x}^{t}$. Notice that the moduli space
$\overline{\mathcal{M}}(\tilde{x},\tilde{x}^{t})$ is only a
pseudo-cycle with boundary in the sense that it is a stratified set
with three strata:
\begin{itemize}
\item[i.]
  a co-dimension two stratum:
  $$
  \Sigma_{\tilde{x},\tilde{x}^{t}}\subset
  \overline{\mathcal{M}}(\tilde{x},\tilde{x}^{t})
  $$
  formed by the bubbled configurations.
\item[ii.]
  a co-dimension one stratum:
  $\partial\mathcal{M}=\cup_{\tilde{z}}\overline{\mathcal{M}}(\tilde{x},\tilde{z})\times
  \overline{\mathcal{M}}(\tilde{z},\tilde{x}^{t})$
\item[iii.]
  a co-dimension zero stratum:
  $\mathcal{M}(\tilde{x},\tilde{x}^{t})~.~$
\end{itemize}

Fix now some $\tilde{x}$ and, to simplify notation, let
$\Sigma=\Sigma_{\tilde{x},\tilde{x}^{t}}$ and notice that $\Sigma$ is
 a compact manifold. The spine map $\sigma$ is defined on
$\overline{\mathcal{M}}(\tilde{x},\tilde{x}^{t})$ with the exception
of $\Sigma$. Notice also that $\Sigma \cap \partial
\mathcal{M}=\emptyset$. Suppose that there exists a continuous
deformation $\sigma'$ of $\sigma$ which agrees with $\sigma$ with the
exception of a neighborhood of $\Sigma$. Then, as
$\overline{\mathcal{M}}(\tilde{x},\tilde{x}^{t})$ is a pseudo-cycle,
the same argument described above for the case
$\tilde{z}\not=\tilde{x}\sharp \alpha_{min}$ applies also here (the
point is that as $\Sigma$ is of co-dimension $2$, the construction of
representing cycles is still possible) and it shows that
$S(\tilde{x},\tilde{x}^{t})$ does not play any role in  $E^{r}$ for
$r\geq 2$. To conclude, our assumption $\dFloer^{\cmin}\circ
\dFloer^{\cmin}\not=0$, implies that there exists at least one
$\tilde{x}$ so that such a deformation $\sigma'$ of $\sigma$ does not
exist. We now want to deduce from this that the first cohomology
group of $\Sigma$ does not vanish.

Given that, by definition $c_{1}(\alpha_{min})=\cmin$, it follows
that each $u\in \Sigma$ is represented by a Floer tube $\R\times
S^{1}$ to which is attached a single sphere in a point
$(t_{u},a_{u})\in \R\times S^{1}$ so that the tube is mapped in $M$
on the constant orbit $\tilde{x}$ and the sphere is mapped to a
pseudo-holomorphic sphere in the class $\alpha_{min}$. Thus, there is
a continuous map
 $$
 \xi : \Sigma\to S^{1}
 $$
so that $\xi(u)=a_{u}$.

To show that $H^{1}(\Sigma ;\Z)\not=0$ it is enough to show that
$\xi$ is not null-homotopic.  Assume that $\xi\simeq 0$. Then $\xi$
can be lifted to an application $\tilde{\xi}:\Sigma\to\R$. Fix
$\chi:\R\to \R$ a smooth function supported on $[-1,1]$ and such that
 $\chi(0)=1$. For $A,s_0$ in $\R$ consider the function
$$
\chi_{s_{0},A}:\R \stackrel{A\chi(s-s_0)}{\longrightarrow} \R\to
S^{1}
$$
where the second map in the composition is $t\to e^{ i t}$.

The graph of this function defines a deformed spine $\Delta (s_0,A) =
 graph (\chi_{s_{0},A})$ on $\R\times S^{1}$ with the property that,
if $A\not= 2k\pi$, it avoids the point $(s_{0},0)$ . For each bubbled
 curve $u\in\Sigma$ we consider the deformed line on the tube given
by  $\Delta_{u}=\Delta (t_{u},\tilde{\xi}(u)+\pi)$. This line avoids
the  point $(t_{u},a_{u})$ and thus avoids the ``bubble''. We obtain
in  this way a continuous spine map: $\sigma':\Sigma\to \Omega M'$
defined by
$$
\sigma'(u)=u(\Delta_{u})~.~
$$
To conclude our proof it is enough to show that this spine map
extends continuously to
$\overline{\mathcal{M}}(\tilde{x},\tilde{x}^{t})$ without modifying
the standard spine map on $\partial \mathcal{M}$. Due to by-now
standard gluing results \cite{Sa}\cite{McSa}, for each point
$x\in\Sigma$ there exists a small neighborhood $U(x)\subset \Sigma$
and a homeomorphism $\phi:\mathbb{C}\times U(x) \to
\overline{\mathcal{M}}(\tilde{x},\tilde{x}^{t})$ so that
$\phi(\{0\}\times U(x))=U(x)$. As $\Sigma$ is compact we can cover it
with a finite number of such neighborhoods which we denote by
$U_{i},\ 1\leq i\leq k$ with corresponding homeomorphisms $\phi_{i}$.
Denote $V_{i}=\phi_{i}(U_{i})$ and let $p_{i}:V_{i}\to U_{i}$ be the
obvious projection. For a point $y\in V(x)$ let
$d_{i}(y)=d(y,p_{i}(y)$ where $d(-,-)$ is (some) distance in
$\overline{\mathcal{M}}(\tilde{x},\tilde{x}^{t})$. By possibly using
smaller neighborhoods $V_{i}$, we may assume that $d_{i}(y)< 1,\
\forall i, y$. Finally, let $h_{i}:U_{i}\to [0,1],\ 1\leq i\leq k$ be
a partition of the unity. We put $U(\Sigma)=\cup V_{i}$. We also
consider a smooth function $\eta:[0,1]\to \R$ which is decreasing,
supported on $[0,1/2]$ and so that $\eta(0)=1$. Let $d'_{i}:V_{i}\to
\R$ be given by $d'_{i}(x)=\eta (d_{i}(x))$.

With these notations we now extend $\sigma'$ to $U(\Sigma)$: we let
$$
 \Delta_{u}=\Delta(\sum_{i}h_{i}(p_{i}(u))t_{p_{i}(u)},\sum_{i} h_{i}(p_{i}(u))d'_{i}(u)(\bar{\xi}(p_{i}(u))+\pi))%
$$
and put $\sigma'(u)=u(\Delta_{u})$. As this map coincides with
$\sigma$ on $\partial U(\Sigma)$ we may extend $\sigma'$ to a
continuous map on all of
$\overline{\mathcal{M}}(\tilde{x},\tilde{x}^{t})$ so that it equals
$\sigma$ outside $U(\Sigma)$. This concludes the proof.

\subsection{Quantum perturbation of the Serre spectral sequence.}
The purpose of this subsection is to show point ii. in Theorem
\ref{theo:spec} and thus conclude the proof of this theorem.

The page $E^{2}$ is well defined and invariant, and by Lemma \ref{lem:e2}
$$
 E^{2}_{p,q}\cong HF_{p}(M; \Lambda)\otimes H_{q}(\Omega (M)).
$$
This is also the first page of the (classical) Serre path-loop spectral
sequence. However, the second differential, $\dFloer^{2}$, on this page
is in general different from the classical one.
To interpret $\dFloer^{2}$ in terms of
binary Gromov-Witten invariants, we will use the construction of
\cite{PSS}. To this end, we start with a quantized-Morse version of the
spectral sequence construction described before.

\subsubsection{The  Quantized - Morse truncated spectral sequence}
\label{subsubsec:Quant-Morse}

To a Morse-Smale pair $(f,g)$ on $M$, together with a generic almost complex
structure $J$ we associate an extended
quantized Morse complex $CM_{\ast}=CM_{*}(f,M, J)$.  This is the free module generated by the
critical points $Crit (f)$ over the ring $\Ring$ together with a differential
which will be described below. The degree of a
critical point $x\in Crit(f)$ is given by its index.

Given the almost complex structure $J$ on $M$, a ``quantum-Morse''
trajectory from a critical point $x$ to a critical point $y$ in class
$\alpha\in\Gamma_{0}$, is a finite collection
$((\gamma_{0},\dots,\gamma_{k}),(S_{1},\dots,S_{k}))$ of paths and
spheres in $M$ such that:
\begin{enumerate}
\item%
each sphere $S_{i}$ is a $J$-holomorphic sphere with a marked real
line $[p_{i,0},p_{i,\infty}]$ on it, and $\sum_{i}[S_{i}]=\alpha$,
\item%
$\forall i$, $\gamma_{i}$ is a piece of flow line of $-\nabla_{g}f$, joining
$S_{i-1}(p_{i-1,\infty})$ to $S_{i}(p_{i,0})$ (with the convention that
$S_{-1}(p_{-1,\infty})=x$ and $S_{k+1}(p_{k+1,0})=y$).
\end{enumerate}
We denote by $\M_{\alpha}(x,y)$ the set of all such objects. For a
generic choice of $(f,g,J)$, it is a smooth manifold of dimension
$$
\dim\M_{\alpha}(x,y)=|x|-|y|+2c_{1}(\alpha)-1~.~
$$
There is no difficulty to prove this as regularity comes down
to the usual transversality of the appropriate evaluation maps
\cite{McSa} (in particular, this is much simpler than the relative
case discussed for example in \cite{BiCo}). To verify the dimension
formula notice that there are not only two marked points on the spheres, but
also a real line joining them.

\medskip

Such a trajectory defines a path from $x$ to $y$ by concatenation of
the flow lines and the marked real lines on the spheres. Notice that
each flow line segment can be parametrized by the value of $-f$,
while on a holomorphic sphere $u:\C\cup\{\infty\}\to\CP^{1}\to M$,
the map $t\in[0,+\infty)\mapsto\int_{|z|\leq t}u^{*}\omega$ is
strictly increasing and defines a parametrization of the marked real
line.

These independent parametrizations of the different segments can now
be shifted and aligned to produce a parametrization of the full
respective path by the segment $[-f(x),-f(y)+\omega(\alpha)]$.
We also assume that the path $\gamma$
used to define $M'$ and turn trajectories into loops
goes through all the critical points of $f$. As a consequence
we obtain a continuous map
$$
\spinemap:\M(x,y)\to\Omega(M').
$$

The space $\M_{\alpha}(x,y)$ of course has a natural Morse-Gromov
compactification $\Mbar_{\alpha}(x,y)$, and the question arises again
of extending $\spinemap$ over it.
Clearly, $\spinemap$ extends continuously over broken trajectories as long
as no bubble components appear (as in \S\ref{plan:Boundary operator}).
However, it might fail to extend over
trajectories where new spheres bubble off. The arguments
used in the discussion of this point for Floer moduli spaces still
apply in this situation, and the map $\spinemap$ can be defined with
the desired continuity and compatibility conditions whenever, as in Lemma \ref{lem:bubbling_cond},
$$|x|-|y|\leq 2\cmin-1~.~$$

Picking chain representatives $m_{\alpha}(x,y)$ of
$\spinemap(\Mbar_{\alpha}(x,y))$, we define a truncated boundary
operator on $CM_{*}$ in the usual way:
$$
  \dFloer x=\sum_{1\leq |x|-|y|+2c_{1}(\alpha)-1\leq 2\cmin-1}
              m_{\alpha}(x,y)\ y e^{\alpha}.
$$
The complex $(CM_{*},\dFloer)$ admits also a differential filtration again defined
by the degree of the elements in $\Z/2<Crit(f)>\otimes\ \Lambda$ and this induces
 a truncated spectral sequence in the same way as before.
 The differential $\dFloer^{2}$ of this spectral sequence has a natural
 interpretation in terms of Gromov-Witten invariants.

\medskip

To see this first notice that if an element $u\in
\Mbar_{\alpha}(x,y)$ with $\alpha\not=0$ is so that it contains $k$
spheres, then the dimension of $\Mbar_{\alpha}(x,y)$ is at least
equal to $k$. Indeed, the choice of the real line on each of the
spheres in $u$ gives rise to a full $S^{1}$ parametric family of
elements in this moduli space. The first consequence of this remark
is that the differential $d^{1}$ in the spectral sequence is simply
$d_{Morse}\otimes id$ which is defined on $\Z/2<Crit(f)>\otimes
H_{\ast}(\Omega M)\otimes \Lambda$ where $d_{Morse}$ is the usual
Morse differential. Indeed, $d^{1}$ involves $0$-dimensional
quantized-Morse moduli spaces and the remark above shows that when
$\alpha\not=0$ these spaces are never $0$-dimensional if non-void.

Suppose now, to shorten the discussion, that $f$ is a perfect Morse
function (if not, critical points should be replaced by a basis of
the Morse homology of $M$).

In this case, for a critical point $x$ of $f$ the second differential $d^{2}x$ is
defined and is given by $d^{2}x=\sum_{\alpha,y}
[m_{\alpha}(x,y)]y$, where the sum is taken over all $(x,y,\alpha)$
such that $\Mbar_{\alpha}(x,y)$ is $1$-dimensional. Notice
that we may associate a homology class in $H_{\ast}(\Omega M)$
to each such moduli space in this case because the Morse differential
vanishes. We let $[m_{\alpha}(x,y)]$ be this class.

In the expression for $\dFloer^{2}$ the sum of the terms
where $\alpha$ is trivial, $d^{2}_{0}$, is given by $1$-parametric
families of Morse trajectories and, as shown in \cite{BaCo1}, this
coincides with the second differential in the Serre spectral sequence
of the fibration $\Omega M\to PM\to M$. On the other hand, when
$\alpha$ is non trivial, as a second consequence of the remark above,
we see that the corresponding moduli space $\Mbar_{\alpha}(x,y)$ is
the set of {\em single} holomorphic spheres in class $\alpha$ with a
marked real line $[p_{0},p_{\infty}]$ such that $p_{0}\in W^{u}(x)$
and $p_{\infty}\in W^{s}(y)$.

The choice of the real line defines an $S^{1}$ action on
$\Mbar_{\alpha}(x,y)$, and the quotients are $0$-dimensional:
$$
 S^{1}\to\Mbar_{\alpha}(x,y)\to\Mbar_{\alpha}(x,y)/S^{1}~.~
$$
Moreover, letting
$$
GW_{\alpha}(x,y)=GW_{\alpha}([x],[y])
$$
be the Gromov-Witten invariants of holomorphic spheres in class
$\alpha$ with two marked points, associated to the homology class of
$x$ and the dual class of $y$, we have:
$$
 GW_{\alpha}(x,y)=\sum\sharp(\Mbar_{\alpha}(x,y)/S^{1}).
$$
In particular, $GW_{\alpha}(x,y)$ is the number of components of
$\Mbar_{\alpha}(x,y)$. Each component of the space $\Mbar_{\alpha}(x,y)$ defines a
loop of loops in $M'$, whose class in $H_{1}(\Omega M')$ is the
image $[\alpha]$ of $\alpha$ under the map $\pi_{2}(M')\to
H_{1}(\Omega M')$. This image is the same for all the components, so
that $[m_{\alpha}(x,y)]= GW_{\alpha}(x,y)[\alpha]\in H_{1}(\Omega M)$. Finally, we
have~:
$$
  d^{2}x=d^{2}_{0}+\sum_{0\not=\alpha\in\pi_{2}(M)}GW_{\alpha}(x,y)[\alpha]ye^{\alpha}.
$$

\subsubsection{Relating the Morse and Floer spectral sequences}

To compare the (truncated) spectral sequences given by the extended
Floer and quantized-Morse complexes, we use the technique introduced
in \cite{PSS} to compare Floer and Morse homologies. For a generic
pair $(H,J)$ we recall the construction of the truncated complex
$C(H,J)$ from \S\ref{plan:Coefficient ring}.

With the notation in that subsection, consider a critical point $x$ of $f$
and a lift $\tilde{y}$ of a (contractible) periodic orbit $y$ of $X_{H}$.

An \emph{hybrid} trajectory from $x$ to $\tilde{y}$ is a quantized-Morse
trajectory - as defined in \S \ref{subsubsec:Quant-Morse} - starting at $x$
but now ending with a disk bounded by $y$.

The definition for a hybrid trajectory are as in in \ref{subsubsec:Quant-Morse},
with the following modifications ~:~
\begin{itemize}
\item[i.]
the last sphere $S_{k}$ is replaced with a disk $u$ with one
cylindrical end, so that in polar coordinates and
away from $0$:
$$
  \C\xrightarrow{u}M \quad\text{with}\quad
  \C=\{0\}\cup\{e^{s+it}, (s,t)\in \R\times S^{1}\}%
$$
\item[ii.]
  the map $u$ satisfies a `cut off' Floer equation. For a fixed cut-off
  function $\chi$ such that $\chi(s)=1$ for $s\geq1$ and $\chi(s)=0$
  for $s\leq0$ we have:
  \begin{equation}\label{eq:Floer cutoff}
    \partial_{s}u+J(u)(\partial_{t}u-\chi(s)X_{H})=0 \ , \ \lim_{s\to+\infty}u(s,t)=y(t)
  \end{equation}
\item[iii.]
  the negative gradient flow arc $\gamma_{k}$ ends at $u(0)$,
\item[iv.]
  the sum of the homotopy class of $u$ with $\sum_{i=1}^{k-1} [S_{i}]$
   defines the capping $\tilde{y}$ of
  $y$.

\end{itemize}


For a generic choice of the data, all the relevant sub-manifolds and
evaluation maps can be made transversal, so that the moduli spaces
$\M(x,\tilde{y})$ of hybrid trajectories are smooth manifolds of
dimension
$$
  \dim\M(x,\tilde{y})=|x|-|\tilde{y}|
$$

These moduli spaces admit a natural compactification - we
refer to \cite{PSS} for the proof. We only recall here that the key
point for showing compactness is to derive a uniform bound
$$
 E(u)=\iint \Vert\frac{\partial u}{\partial s}\Vert^{2}dsdt
 \leq \action(\tilde{y})+\Vert H\Vert_{\infty}.
$$
for the energy from the ``cut off'' Floer equation.

To turn a hybrid trajectory in $\M(x,\tilde{y})$ into a path  from
$x$ to $\tilde{y}(0)$, it is enough to choose a parametrization of
the real line $u(\R)$ on the terminal disk and for that we may use
the energy of the curve:
\begin{equation}
  \label{eq:energy}
E(r)=\int_{(-\infty,r]\times S^{1}}\Vert\frac{\partial u}{\partial s}\Vert^{2}dsdt~. %
\end{equation}

This choice defines a continuous spine map
$\spinemap:\M(x,\tilde{y})\to\Omega(M')$, that can again be extended
to the natural compactification $\Mbar(x,\tilde{y})$ up to dimension
$2\cmin-1$. Picking compatible chain representatives of these spaces,
we obtain chains $m(x,\tilde{y})\in C_{*}(\Omega(M'))$ such that:
\begin{multline}
\dChain m(x,\tilde{y})=\sum_{\substack{
      0\leq |x|          -|ze^{\alpha}|-1\leq 2\cmin-1\\
      0\leq |ze^{\alpha}|-|\tilde{y}|    \leq 2\cmin-1}}
          m(x,ze^{\alpha})m(ze^{\alpha},\tilde{y})
+\\+
 \sum_{\substack{
      0\leq |x|        -|\tilde{z}|  \leq 2\cmin-1\\
      0\leq |\tilde{z}|-|\tilde{y}|-1\leq 2\cmin-1}}
          m(x,\tilde{z})m(\tilde{z},\tilde{y})
\end{multline}

(where $m(x,ze^{\alpha})=m_{\alpha}(x,z)$). Consider now the truncated morphism
$\phi$ given by
$$
\phi(x)=\sum_{|x|-|\tilde{y}|\leq 2\cmin-1} m(x,\tilde{y})\tilde{y}~.~
$$
As expected, the map $\dFloer\phi-\phi\dFloer$ fails to vanish in
general, but one easily checks that $(\dFloer\phi-\phi\dFloer)(x)$ is
supported on elements $\tilde{y}$  with
$|\tilde{y}|\leq |x|-2\cmin$. This means that $\phi$ induces a
morphism $\Phi$ between the respective truncated spectral sequences
of order $c_{min}$.

\

Notice that the $\Phi^{1}$ coincides with:
$$\phi'\otimes id:C_{Morse}(f,g)\otimes H_{\ast}(\Omega M)\to CF_{\ast}(H)\otimes H_{\ast}(\Omega M)$$
where $\phi'$ is the usual PSS morphism and $C_{Morse}(f,g)$ is the
Morse complex of $(f,g)$. But, as $\phi'$ induces an isomorphism in
homology, this implies that $\Phi^{2}$ is an isomorphism which, in
particular, proves the point ii. of Theorem \ref{theo:spec} and
concludes the proof of this theorem.

\section{Examples, applications and further comments}
\subsection{Extensions}\label{subsec:ext} We recall that the setting considered till now in the paper
was that of a closed, simply-connected, monotone manifold for which
$c_{min}\geq 2$. All the constructions described previously in the
paper extend much beyond this setting. We will discuss here a few
such generalizations.

\subsubsection{$\pi_{1}\not=0$.}
There are two essential ways to perform our constructions in the presence of a non trivial  fundamental
group. They both stem from the fact that the only place where the fundamental group of $M$  affects
the construction is in the possible dependence of the resulting homology on the path
$\gamma$ which is used to define  the quotient $$M\to M'$$ as described
in \S\ref{plan:Coefficient ring}. Of course, at the level of the spectral sequences
$\pi_{1}(M)\not=0$ also plays a role as local coefficients might be necessary.

A. The first way to deal with the fundamental group consists in
enlarging the Novikov ring by tensoring with
the group ring $\Z/2\ [\pi_{1}(M)]$. Geometrically, this can be viewed as performing
all the topological constructions on the universal covering, $\widetilde{M}$, of $M$
even though all equations satisfied by the elements in our new moduli spaces take place
after projection into $M$. The covering $\CovOrbits_{H}$ is
replaced by the covering $\CovOrbits_{H}'$ which is the pull-back of
$\widetilde{M}\to M$ over $\CovOrbits_{H}\to \Orbits_{H}\to M$.
In this case our truncated complex is isomorphic to:
$$\Z/2<\Orbits_{H}>\otimes \Lambda\otimes \Z/2[\pi_{1}(M)]\otimes C_{\ast}(\Omega M)~.~$$

\

B. A second possibility is {\em localization} or change of coefficients. This is maybe even more
useful in applications than A and consists in replacing in all the construction the
coefficient ring $C_{\ast}(\Omega M)$ by $C_{\ast}(\Omega X)$ where $X$ is some simply-connected
 topological space which is endowed with a map:
 $$\eta:M\to X~.~$$
All our moduli spaces are represented inside $\Omega(M')$ and, by composition
with the map $\Omega \eta:\Omega(M')\to \Omega(X)$, they are also represented
inside $\Omega (X)$. The results in Theorem \ref{theo:spec} remain true after this
change of coefficients except that $H_\ast (\Omega M)$ is replaced by
$H_{\ast}(\Omega X)$ and the path loop fibration over $M$ is replaced with the
fibration of base $M$ which is obtained by pull-back over the map $\eta$ from the
path-loop fibration over $X$, $\Omega X\to PX\to X$.

\subsubsection{$c_{min}= 1$.}\label{subsubsec:c1} It is easy to see that even if $c_{min}=1$ the $E^{2}$
term of our spectral sequence is well defined together with the map $d^{2}$ (which might not
be a differential though) and Theorem \ref{theo:spec} remains true
for the $E^{2}$ term. This happens because to prove the invariance of $d^{2}$ only
moduli spaces of dimension $2$ are needed. In turn, as bubbling is a codimension two
phenomenon this means that the bubbling points can be avoided when defining the
spine map over these moduli spaces (as also discussed in Remark \ref{rem:non-monotone}).
\subsubsection{Rational coefficients.} The moduli spaces we use in this paper admit coherent
orientations and by taking these into account we may replace everywhere $\Z/2$ with $\Q$.
\subsubsection{Lack of monotonicity.} It is expected that Theorem \ref{theo:spec}
remains true for the $E^{2}$ term of the spectral sequence even if $(M,\omega)$ is not monotone
(see also Remark \ref{rem:non-monotone}).
Of course, in this case multi-valued perturbations are needed and, thus, the use of rational
coefficients is mandatory.
\subsubsection{Non-compactness.} Finally, it is obviously possible to extend this
theory to the case when $M$ is not compact if it is convex at
infinity. In that case the Hamiltonians used should have compact
support.

\subsection{Examples}
\subsubsection{$\CP^{1}$}
Take now $M=\CP^{1}$, and consider the Morse function having only one
maximum $a=\infty$ and one minimum $b=0$ as critical points. Let us
pick a simple path from $b$ to $a$ to serve as the base point of
$M'$. As auxiliary data, we can simply stick to the standard metric
and complex structure on $\CP^{1}$, and use no perturbations at all~:
one easily checks that genericity is fulfilled for all the moduli
spaces involved in the computations below.

The page $2$ of the spectral sequence is simply
$H_{*}(\CP^{1})\otimes H_{*}(\Omega
S^{2})$. Let $\alpha$ denote the identity map $S^{2}\to\CP^{1}$. Seen
as the $S^{1}$ family of flow lines going from $a$ down to $b$,
$\alpha$ defines a cycle $[\alpha]$ that generates
$H_{1}(\Omega(S^{2}))$.

The Novikov ring is generated by the multiples of $\alpha$~:
$$
  \Lambda=\left\{             \sum_{\lambda_{k}\in\Z_{2}
  }\lambda_{k}e^{k\alpha},\           \right\}.
$$
and we have $c_{1}(\alpha)=2$.

To make the differential more explicit, we will ``unfold'' the
spectral sequence by removing the Novikov ring from the coefficients,
and thinking of $\{ae^{k\alpha}\}$ or $\{be^{k\alpha}\}$ as free families.

To compute the differential $\dFloer_{2}$, we have to compute all the
$1$-dimensional moduli spaces. Because of the invariance of the
moduli spaces under the action of $\pi_{2}(S^{2})$ on both ends of
the trajectories, we can restrict to spaces of the form
$\Mbar(x,ye^{k\alpha})$ with $x,y\in\{a,b\}$. The dimension of this
space is
$$
\dim \Mbar(x,ye^{k\alpha})=|x|-|y|+4k-1
$$
so there are only two possibilities~:
\begin{itemize}
\item[-]
  $k=0$, $x=a$ and $y=b$,

\item[-]
  $k=1$, $x=b$ and $y=a$.
\end{itemize}

The first moduli space consists in classical flow lines only~: it
contributes to the classical part $\dFloer^{2}_{0}$ of $\dFloer^{2}$,
and we have~:
$$
\dFloer_{0}^{2}(a)=[\alpha]b, \quad \dFloer_{0}^{2}(b)=0,
$$
so that the page $2$ of the ``classical'' spectral sequence (tensored
by the Novikov ring) has the following form~:
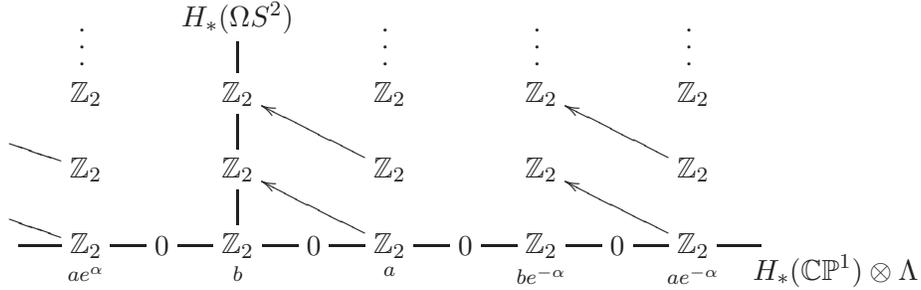
$$
\xymatrix@!C=4ex@!R=1ex{
&&&{H_{*}(\Omega S^{2})}\ar@{-}[d]&&&&&&&\\
& \Z_{2}\ar@{{}*{\cdot}{}}[u]&& \Z_{2}\ar@{-}[d]           &&
\Z_{2}\ar@{{}*{\cdot}{}}[u]&& \Z_{2}\ar@{{}*{\cdot}{}}[u]&&
\Z_{2}\ar@{{}*{\cdot}{}}[u]&
\\
& \Z_{2}\ar@{-}[ul]+<0ex,-4ex>&& \Z_{2}\ar@{-}[d]            &&
\Z_{2}\ar[ull]&& \Z_{2}&& \Z_{2}\ar[ull]&
\\
  \ar@{-}[r]&
  \Z_{2}\ar@{-}[ul]+<0ex,-4ex>\ar@{-}[r]\ar@{{}{}{}}@<-1ex>_{ae^{
  \alpha}}[]&0\ar@{-}[r]& \Z_{2}\ar@{-}[r]\ar@{{}{}{}}@<-1ex>_{b }[]
  &0\ar@{-}[r]& \Z_{2}\ar@{-}[r]\ar@{{}{}{}}@<-1ex>_{a
  }[]\ar[ull]&0\ar@{-}[r]&
  \Z_{2}\ar@{-}[r]\ar@{{}{}{}}@<-1ex>_{be^{-\alpha}}[] &0\ar@{-}[r]&
  \Z_{2}\ar@{-}[r]\ar@{{}{}{}}@<-1ex>_{ae^{-\alpha}}[]\ar[ull]&
  \save[]+<0ex,-2ex>*[r]{H_{*}(\CP^{1})\otimes\Lambda}\restore }
$$
\medskip

The second moduli space, $\Mbar(b,ae^{\alpha})$, involves holomorphic
spheres of degree $1$, and determines the quantum component
$\dFloer^{2}_{Q}$ of $\dFloer^{2}$. Since there are no flow lines
going out of $b$ or into $a$, it consists in holomorphic spheres of
degree $1$ with a marked real line from $b$ to $a$. This is the same
cycle as $\alpha$, but with reversed orientation. as a consequence,
we have
$$
\dFloer^{2}_{Q}(a)=0\text{ and } \dFloer^{2}_{Q}(b)=ae^{\alpha},
$$
and the page $2$ of the full spectral sequence has the following
form~:
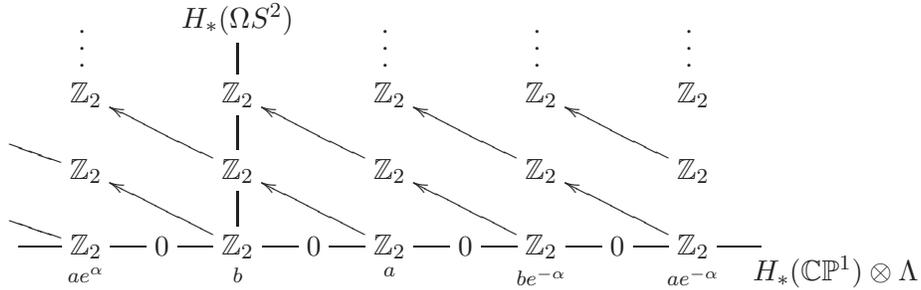
$$
\xymatrix@!C=4ex@!R=1ex{
&&&{H_{*}(\Omega S^{2})}\ar@{-}[d]&&&&&&\\%
  &
  \Z_{2}\ar@{{}*{\cdot}{}}[u]&&
  \Z_{2}\ar@{-}[u]\ar@{-}[d]&& \Z_{2}\ar@{{}*{\cdot}{}}[u]&&
  \Z_{2}\ar@{{}*{\cdot}{}}[u]&& \Z_{2}\ar@{{}*{\cdot}{}}[u]&%
\\%
  & \Z_{2}\ar@{-}[ul]+<0ex,-4ex>&& \Z_{2}\ar[ull]\ar@{-}[d]&& %
 \Z_{2}\ar[ull]&& \Z_{2}\ar[ull]&& \Z_{2}\ar[ull]&%
\\
  \ar@{-}[r]& \Z_{2}\ar@{-}[ul]+<0ex,-4ex> %
  \ar@{-}[r]\ar@{{}{}{}}@<-1ex>_{ae^{ \alpha}}[]        &0\ar@{-}[r]&%
  \Z_{2}\ar@{-}[r]\ar@{{}{}{}}@<-1ex>_{b }[]\ar[ull]&0\ar@{-}[r]& \Z_{2}\ar@{-}[r]\ar@{{}{}{}}@<-1ex>_{a }[]\ar[ull]&0\ar@{-}[r]& \Z_{2}\ar@{-}[r]\ar@{{}{}{}}@<-1ex>_{be^{-\alpha}}[]\ar[ull]&0\ar@{-}[r]& \Z_{2}\ar@{-}[r]\ar@{{}{}{}}@<-1ex>_{ae^{-\alpha}}[]\ar[ull]& \save[]+<0ex,-2ex>*[r]{H_{*}(\CP^{1})\otimes\Lambda}\restore }
$$

Notice that $(\dFloer^{2})^{2}a=[\alpha^{2}]ae^{\alpha}\neq 0$. So
some bubbling has to occur on a $3$ dimensional moduli space. And in
fact, the moduli space $\Mbar(a,ae^{\alpha})$ is $3$ dimensional and
consists in flow lines going out of $a$ down to some point $p$, and a
holomorphic sphere of degree $1$ with a marked real line from $p$ to
$a$. When the point $p$ goes to $b$, the flow line brakes, and we see
that this space is involved in the computation of
$(\dFloer^{2})^{2}(a)$. On the other hand, when the point $p$ goes to
$a$, we are left with the constant trajectory from $a$ to itself,
with an (unparametrized) holomorphic sphere attached to it. Here, the
critical point $a$ is seen as a constant tube, with a marked real
line~: this marked line is responsible for the bubbling monodromy.

It is interesting to note that this bubbling is in fact equivalent
to the fact that in the Pontryaguin algebra $H_{\ast}(\Omega S^{2};\Z/2)$
the non-vanshing class in $H_{1}(\Omega S^{2};\Z/2)$ has a
non-vanishing square.

%

\subsubsection{$\CP^{n}$ for $n>1$}
The computation can be achieved in the same way on $\CP^{n}$ for
$n>1$. Notice that the minimal first Chern class is $n+1\geq 3$ so
that the spectral sequence still exists after page $2$. It is an easy
computation to see that the quantum component of the differential
$d^{2}$ is given by~:
$$
d^{2}_{Q}[pt]=[\Delta]\otimes[\CP^{n}]e^{\Delta}
$$
where $\Delta$ is a complex line in $\CP^{n}$.
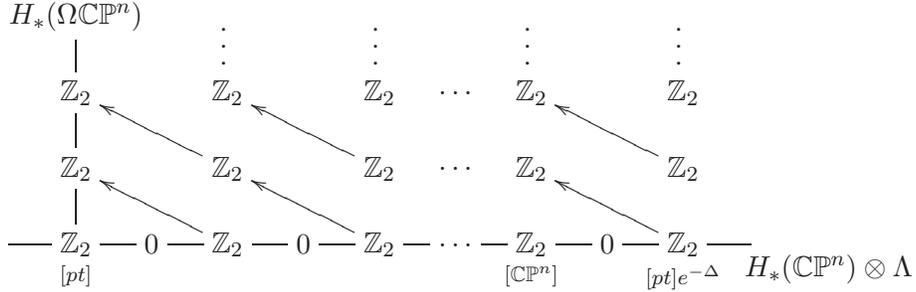
$$
\xymatrix@!C=4ex@!R=1ex{
&{H_{*}(\Omega \CP^{n})}\ar@{-}[d]&&&&&&&&\\
  & \Z_{2}\ar@{-}[u]\ar@{-}[d]&& \Z_{2}\ar@{{}*{\cdot}{}}[u]&&
  \Z_{2}\ar@{{}*{\cdot}{}}[u]&\dots & \Z_{2}\ar@{{}*{\cdot}{}}[u]&&
  \Z_{2}\ar@{{}*{\cdot}{}}[u]&
\\
  & \Z_{2}\ar@{-}[d]&& \Z_{2}\ar[ull]&& \Z_{2}\ar[ull]&\dots &
  \Z_{2}&& \Z_{2}\ar[ull]&
\\
  \ar@{-}[r]& \Z_{2}\ar@{-}[r]\ar@{{}{}{}}@<-1ex>_{[pt] }[]&
  0\ar@{-}[r]& \Z_{2}\ar@{-}[r]\ar@{{}{}{}}@<-1ex>_{ }[]\ar[ull]&
  0\ar@{-}[r]& \Z_{2}\ar@{-}[r]\ar@{{}{}{}}@<-1ex>_{ }[]\ar[ull]&
  \dots\ar@{-}[r]& \Z_{2}\ar@{-}[r]\ar@{{}{}{}}@<-1ex>_{[\CP^{n}]
  }[]& 0\ar@{-}[r]&
  \Z_{2}\ar@{-}[r]\ar@{{}{}{}}@<-1ex>_{[pt]e^{-\Delta}}[]\ar[ull]&
  \save[]+<0ex,-2ex>*[r]{H_{*}(\CP^{n})\otimes\Lambda}\restore }
$$
In particular, the pages of the spectral sequence all vanish after
page $2$.

The contrast between the situation $n=1$ and $n>1$ comes
from the properties of the Pontryaguin product in
$H_{1}(\Omega\CP^{n})$. This product is involved in the computation
of ${\dFloer^{2}}\circ \dFloer^{2}$, in particular~:
$$
\dFloer^{2}(\dFloer^{2}([pt]))=[\Delta]*[\Delta]\otimes a_{n-1}
e^{\Delta}
$$
where $a_{n-1}$ is a generator of $H_{2(n-1)}(\CP^{n})$. What is
truly remarkable here is that as $c_{min}=n+1$, when $n> 1$,
our construction of the truncated spectral sequence shows
that $\dFloer^{2}\circ\ \dFloer^{2}=0$ which implies
$[\Delta]\ast[\Delta]=0$ in the Pontryaguin ring. Of course, this
relation is well-known by purely topological methods but it is remarkable
that it is a consequence of the existence of the quantized Serre
spectral sequence.
Moreover, by Theorem \ref{theo:spec} ii,
$\dFloer^{2}$ can be
expressed in terms of Gromov-Witten invariants together with classical
Hopf ones and this discussion shows that the relations among them
in the Pontryaguin algebra are not trivial.

\subsection{Fibrations over $S^{2}$}

Given a loop $\phi$ in $\Ham(M)$, one can construct a fibration
$E_{\phi}$ over $S^{2}$, obtained by gluing two trivial fibrations
over the disk via $\phi$.

P.~Seidel \cite{Seidel} used sections of this fibration to associate  an
invertible endomorphism on $H_{*}(M)$ to each such $\phi$, deriving
strong topological restrictions on elements in $\pi_{1}(\Ham(M))$.
We first give an outline of the construction of this morphism in the
context of Morse homology and then explain how it is related to our
construction .

Let $\Omega$ be a symplectic form on $E_{\phi}$ such that its
restriction to the fibers is (cohomologuous to) $\omega$, and let
$J_{\phi}$ be an almost complex structure $\Omega$-compatible on
$E_{\phi}$.

Let $f:M\to\R$ be a Morse function on $M$, and let $\tilde{f}$ be a
Morse function on $E_{\phi}$ such that
\begin{itemize}
\item[-]
$\tilde{f}(z,m)=f_{M}(m)+|z^{2}|+cst$  in a local chart near $0$
\item[-]
$\tilde{f}(\tau,m)=f_{M}-|\tau|^{2}$ in a local chart near $\infty$
\item[-]
$\tilde{f}$ has no other critical point than those in the fibers of $0$ and
$\infty$.
\end{itemize}

If $x$ is a critical point of $f$, we denote by $x_{+}$ and $x_{-}$
the corresponding critical points above $\infty$ and $0$
respectively.

We have $i(x_{+})=i(x)+2$ and $i(x_{-})=i(x)$.

For the purpose of our discussion, notice that $J_{\phi}$ can be
chosen to be the product
$\big(\begin{smallmatrix}i&0\\0&J\end{smallmatrix}\big)$
of structures on $S^{2}$ and $M$ in local charts $U_{0}$ and
$U_{\infty}$ near the fibers over $0$ and $\infty$: the only curves
contained in this region are in fact contained in a fiber, and the
almost complex structure is regular for them. The almost complex
structure $J_{\phi}$ can then be made regular for all curves by
perturbing it in the complement of this region.

Roughly speaking, the Seidel morphism is obtained by considering
$0$-dimensional moduli spaces of flow lines going out of a critical
point $x_{-}$, hitting a $J_{\phi}$ holomorphic section of
$E_{\phi}$, followed by a second flow line flowing from the section
down to a critical point $y_{+}$.

To be able to compare homology classes of section with homology
classes in $M$ we fix a section $s_{0}$ of the fibration $E_{\phi}$:
the homology classes having degree $1$ over the base are then the
classes of the form $s_{0}+i_{*}\alpha$, for $\alpha\in H_{2}(M)$
where $i:H_{2}(M)\to H_{2}(E_{\phi})$ is induced by the inclusion
$M=M\times\{0\}\hookrightarrow E_{\phi}$. The Seidel morphism is
induced at the homology level by the map $\Phi$:
$$
  \Phi(x)=\sum_{i(x)-i(y)+2c_{1}(s_{0}+i_{*}\alpha)-2=0}
  GW(x_{-},y_{+};s_{0}+i_{*}\alpha)ye^{\alpha}
$$

We now discuss how to interpret this morphism as a component of the
differential $\dFloer^{2}$ the truncated spectral sequence associated
to $E_{\phi}$. The version we will use is a variant of the quantized
Morse one from \S\ref{subsubsec:Quant-Morse}. With the notation in
\S\ref{subsubsec:Quant-Morse}, we  write the differential of the
quantized Morse complex $CM(\tilde{f},E_{\phi})$ as $\dFloer
x=\sum_{k} d_{k}(x)$ where
$$d_{k}=\sum_{\lambda, \ deg(\lambda)= k }m_{\lambda}(x,y)y e^{\lambda}$$
with the degree considered over the base. This decomposition induces
an analogue one for the differentials of the associated truncated
spectral sequence which we will denote by
$$d^{r}=\sum d^{r;k}$$ with $d^{r;k}$ induced by $d_{k}$.

For $k=0$, all the moduli spaces involved in $\dFloer^{2;0}x_{-}$ lie
in the same fiber as $x_{-}$: they are all images of the
corresponding moduli spaces in $M$ via the inclusion $i$ of $M$ in
$E_{\phi}$ as the fiber over $0$. At the homology level, we have the
following commutative diagram:
$$
\xymatrix{%
H_{*}(M)\ar^{i_{*}}[d]\ar^>(.6){\dFloer^{2}}[rr]&&
H_{*}(M)\otimes H_{1}(\Omega(M))\ar^{i_{*}}[d]\\
H_{*}(E_{\phi})\ar^>(.6){\dFloer^{2;0}}[rr]&& H_{*}(E_{\phi})\otimes
H_{1}(\Omega(E_{\phi}))
}%
$$

\medskip

Consider now the case $k=1$. For dimensional reasons, $1$-dimensional
moduli spaces of degree $1$ quantum trajectories starting at a point
$x_{-}$ have to end in a point $y_{+}$.

By Theorem \ref{theo:spec}, the differential $\dFloer^{2;1}$ applied to the critical point
$x_{-}$ has the following form:
\begin{gather*}
\dFloer^{2;1} x_{-}=\sum_{\alpha\in H_{2}(M)} GW(x_{-},y_{+};s_{0}+i_{*}\alpha)\ [s_{0}+i_{*}\alpha]\ y_{+} e^{s_{0}+i_{*}\alpha}.
\end{gather*}

Using $\pi:E_{\phi}\to S^{2}$ to change coefficients and replace
$\Omega(E_{\phi})$ by $\Omega(S^{2})$, and observing that the classes
$[s_{0}+i_{*}\alpha]$ are all sent to the generator $\alpha$ of
$H_{1}(\Omega(S^{2}))$, we get the following commutative diagram:
$$
\xymatrix{%
H_{*}(M)\ar^{i_{*}}[d]\ar^{\Phi}[rr]&&
H_{*}(M)\ar^>(.7){\mathop{\mathrm{Id}}\otimes[\alpha]}@{^{(}->}[r]&
H_{*}(M)\otimes H_{1}(\Omega(S^{2}))
\\
H_{*}(E_{\phi})\ar^>(.7){\dFloer^{2;1}}[rr]&& H_{*}(E_{\phi})\otimes H_{1}(\Omega(E_{\phi}))%
\ar_{\cap[M]\otimes \pi}[ur]
}%
$$
This relates the Seidel morphism $\Phi$ and the $\dFloer^{2;1}$
component of the differential of the spectral sequence. From this
point of view, when they exist, the higher dimensional components
$\dFloer^{r;1}$ can be viewed as higher dimensional analogues of the
Seidel morphism.

\subsection{Non trivial periodic orbits for Morse functions}

The construction of the truncated spectral sequence can be used
to exhibit extra periodic orbits for Morse functions in some
particular situations.

\subsubsection{Proof of Corollary \ref{thm:extra periodic orbits}}

Let $(M,\omega)$ be a monotone symplectic manifold, and consider a perfect
Morse function $f$ on $M$ which is self indexed. From the statement of the Corollary
recall  that there are two Morse homology classes $x$, $z$, $|z|>|x|$,
which are $\dFloer^{r}$-related.
Due to the self-indexing condition each of these classes is represented
by a linear combination of critical points with the same
critical value, $x=\sum_{i}x_{i},\ z=\sum_{i}z_{i}$. We let $f(x)=f(x_{i})$.
Recall also
that we assume that $H_{k}(M)\otimes \Lambda_{q}=0$ for $|z e^{\lambda}|< k+q < |x|$.

For $A\in\R$ large enough, consider a smooth increasing function
$\phi_{A}:\R\to\R$, such that:
\begin{itemize}
\item[-] $\phi_{A}(t)=t$ for $ t\leq f(x)+1/(2A)$
\item[-] $\phi_{A}(t)=t+A$ for $ t\geq f(x)+1/A$
\end{itemize}
For any $A$, the function $f_A=\phi_{A}\circ f$ has the same
critical points as $f$ with the same (unparametrized) flow lines, but
 the critical levels above $f(x)$ are shifted upward. The critical points
of the same index continue to share the same level hypersurface.
Of course, the existence of non-trivial characteristics for $f$ and $f_{A}$
is equivalent.

Suppose that $f_A$ has no non trivial periodic orbit. It is then
easy to see that, by possibly composing $f$ with another diffeomorphism
$\R\to \R$ whose effect is to diminish the size of the derivatives of
$f$ close to its critical values, and using
a generic almost complex structure which is time-dependent, Floer
theory may be applied to the Hamiltonian $f_{A}$ and the
Conley-Zehnder index of the critical points agrees with their Morse
index. Moreover, the (not extended) Floer and Morse complexes are
then the same (indeed, as the homology of the Floer complex has to be
isomorphic with Morse homology it follows in this case that the Floer
differential is also trivial). Thus, $x$ and $z$ also give Floer
homology classes. We now pick the constant $A$ so that $A\geq \rho
(2n+r)$ where $\rho$ is the monotonicity constant ($\omega
(\alpha)=\rho c_{1}(\alpha)$).

Consider the truncated spectral sequence associated to the
Hamiltonian $X_{f_{A}}$. By hypothesis we know that $x$ is $\dFloer^{r}$-
related to $ze^{\lambda}$. This implies that there is a critical point
$z_{j}$ so that there are Floer trajectories
from one of the $x_{i}$'s to $z_{j}e^{\lambda}$. Indeed, as
$|x|> k+q > |ze^{\lambda}|$ implies $H_{k}(M)\otimes \Lambda_{q}=0$ the differential
$\dFloer^{r}$ is the first one relating the vertical line
through $|x|$ to the one through $|ze^{\lambda}|$. In other words,
letting $p=|ze^{\lambda}|$, we have that
$E^{r}_{p,\ast}$ is a subgroup of $E^{2}_{p,\ast}$ and so, if there are
 no flow lines relating some $x_{i}$ to a $z_{j}e^{\lambda}$, then $x$ and $ze^{\lambda}$ can not be
$\dFloer^{r}$-related.

In view of this we have:
$$|x|-r = |z e^{\lambda}|=|z|-2c_{1}(\lambda)$$
which means that $c_{1}(\lambda)\leq n +r/2$ and  also
$$f(x_{i})=f_{A}(x_{i})\geq f_{A}(z_{j})-\omega(\lambda)=f(z_{j})+A-\omega(\lambda)~.~$$

But given of our choice of constant $A$,
$$f(z_{j})+A-\omega(\lambda)\geq f(z_{j})+A - \rho(2n+r) > f(x_{i})$$
which leads to a contradiction and concludes the proof.

\qed

\

The same technique applies in many other variants of the
  situation described above. The basic idea is to ensure the existence of
  a sequence of trajectories, ``ending at a higher level than its
  starting point'' (there was just one such trajectory in the case above)
  in such a way that the relevant intermediate points can be shifted out
  of the action window (as done before using $\phi_{A}$).
For this, besides identifying a chain of differentials which relate a
succession of homology classes in the spectral sequence one also needs to be able
to choose appropriate chains representing these classes (the self indexing
condition and the homological ``gap" condition had this purpose above).
Here is such a variant valid when
$(\dFloer^{\cmin})^{2}\neq 0$.

\begin{corollary}\label{claim:case dcmin^2!=0}
  If $(\dFloer^{\cmin})^{2}\neq 0$, any self-indexed Morse function on $M$ has
  at least one closed characteristic.
\end{corollary}

\begin{proof} Fix  a self-indexed Morse function $f:M\to \R$. We assume
that it has no non-trivial close characteristics.

Let $\xi \in E^{c_{min}}_{p,q}$ be a class so that
$(d^{c_{min}})^{2}(\xi)\not=0$.

As in the previous proof we may assume that the critical points
of $f$ are non-degenerate periodic orbits of $X_{f}$ and their
Conley-Zehnder index agrees with the Morse index. Thus we may apply
our construction of the truncated quantized Serre spectral sequence
to the Hamiltonian $f$ (together with a generic time-dependent almost
complex structure). We may also assume, after possibly composing $f$
with an appropriate diffeomorphism $\R\to \R$, that:

\begin{itemize}
\item[$\ast$]for any critical points $x\in Crit(f)$
the interval $[f(x),f(x)+\rho(n)+\omega_{min}]$
does not contain any critical values different from $f(x)$.
\end{itemize}

Here $\rho$ is as before the monotonicity constant so that
$\omega_{min}=\rho c_{min}$.
From the discussion in \S\ref{subsec:monodromy} we see that
$(d^{c_{min}})^{2}(\xi)\not=0$ implies that for some critical point $x\in Crit (f)$
we have that the moduli space $\overline{\mathcal{M}}(x,x\# \alpha_{min})$ is non
void and has a non-void codimension one stratum $\Sigma_{1}$ consisting of broken Floer trajectories
as well as a non-void codimension two stratum $\Sigma_{2}$
consisting of Floer trajectories with some bubble attached.

Assume that among the broken trajectories in $\Sigma_{1}$ there is
one which joins $x$ to $ye^{\alpha}$ followed by a second trajectory
from $ye^{\alpha}$ to $xe^{\alpha_{min}}$.

We then have:
$$|x|-|y|+2c_{1}(\alpha)-1\geq 0, \ |y|-|x|+2c_{1}(\alpha_{min}-\alpha)-1\geq 0~.~$$

Notice that this implies that $|y|\not= |x|$. Indeed, if $|y|=|x|$,
the first inequality implies that $c_{1}(\alpha)>0$ and the second that $c_{1}(\alpha)< c_{min}$
which is not possible.

There is also an inequality involving the actions:
\begin{equation}\label{eq:action}
f(x)\geq f(y)-\omega(\alpha)\geq f(x)-\omega_{min}~.~
\end{equation}
There are two cases to consider now. If $|y|>|x|$, then $c_{1}(\alpha)> 0$ so that
$c_{1}(\alpha)\geq c_{min}$
and we also need to have $2 c_{1}(\alpha)\leq |y|-|x|-1 +2c_{min}< 2n + 2c_{min}$.
By monotonicity this means $f(y)-\omega(\alpha)> f(y)-\rho(n+c_{min})$.
Recall that $f$ is self-indexed, $|y|> |x|$ as well as our assumption $\ast$ on the
critical values of $f$. This implies that $f(y)\geq f(x)+\rho(n+c_{min})$
which contradicts the first inequality in (\ref{eq:action}). The second
case is $|y|< |x|$.  Then $c_{min}\geq c_{1}(\alpha)> -2n$.
This means $f(y)-\omega(\alpha)\leq f(y)-\omega_{min}< f(x)-\omega_{min}$
which contradicts the second inequality in (\ref{eq:action}) and
concludes the proof.

\end{proof}

\end{document}